\documentclass[]{article}
\usepackage{amsmath,amsthm}
\usepackage{amssymb}
\usepackage{fullpage}
\usepackage{enumitem}
\usepackage{tikz}
\usepackage{hyperref}
\hypersetup{
colorlinks   = true,    
urlcolor     = blue,    
linkcolor    = blue,    
citecolor    = red      
}

\newcommand{\ie}{{\it i.e.}}

\newcommand{\ones}{\mathbf 1}
\newcommand{\zeroes}{\mathbf 0}
\newcommand{\reals}{{\mbox{\bf R}}}

\newcommand{\diag}{\mathop{\bf diag}}

\newcommand{\Expect}{{\mathbf E}}

\newcommand{\indep}{\perp \!\!\! \perp}

\newcommand{\paren}[1]{\ensuremath{\left(#1\right)}}
\newcommand{\bmat}[1]{\ensuremath{\begin{bmatrix}#1\end{bmatrix}}}
\newcommand{\pmat}[1]{\ensuremath{\begin{pmatrix}#1\end{pmatrix}}}

\title{Strategic Asset Allocation with Illiquid Alternatives}
\author{Eric Luxenberg  \and Stephen Boyd \and Mykel Kochenderfer \and Misha van Beek \and Wen Cao  \and Steven Diamond \and Alex Ulitsky \and Kunal Menda \and Vidy Vairavamurthy}
\begin{document}

\maketitle
\begin{abstract} We address the problem of strategic asset
allocation (SAA) with portfolios that include illiquid alternative asset 
classes.
The main challenge in portfolio construction with illiquid asset 
classes is that we do not have direct control over our positions,
as we do in liquid asset classes.  Instead we can only make
commitments; the position builds up over time as capital calls come in,
and reduces over time as distributions occur, neither of which the investor has
direct control over.
The effect on positions of our commitments
is subject to a delay, typically of a few years, and is also 
unknown or stochastic.
A further challenge is the requirement that we can meet the capital calls,
with very high probability, with our liquid assets.

We formulate the illiquid dynamics as a random linear system,
and propose a convex optimization based model predictive control (MPC)
policy for allocating liquid assets and making new illiquid commitments
in each period.
Despite the challenges of time delay and uncertainty,
we show that this policy attains performance surprisingly close to 
a fictional setting where we pretend the illiquid asset classes
are completely liquid, and we can arbitrarily and immediately
adjust our positions.
In this paper we focus on the growth problem, with no external 
liabilities or income, but the method is readily extended to handle 
this case.
\end{abstract}

\clearpage
\tableofcontents
\clearpage

\section{Introduction} 
There is considerable investor interest across several
financial contexts in constructing portfolios which mix liquid and illiquid
assets, especially illiquid alternative investments. 
We wish to perform strategic asset allocation to asset classes that 
include illiquid alternative assets, as well as more liquid asset classes.
Several challenges arise.  First, we can only augment our illiquid
positions by making capital commitments. Moreover, these commitments only indirectly affect
our illiquid position through uncertain and delayed capital calls, that we have no
direct control over.
A further challenge is the solvency requirement: we should be able to fund
the capital calls from our liquid positions with very high probability.
A simple strategy to guarantee coverage of capital calls
is to keep an amount equal to the uncalled capital commitments in cash.
However this creates significant cash drag,
since this cash could be invested in higher returning liquid assets.
The method we describe in this paper addresses all of these issues.

\section{Previous work}
There is a rich history of studying portfolio construction. Our work helps
extend the modern portfolio theory framework developed by
Markowitz~\cite{Markowitz1952} and Merton~\cite{Merton1969}, which focuses on
liquid assets. We contribute to the further study of illiquidity and
multi-period planning. While this work takes as an input a stochastic model
which describes the risk and return of illiquid investments, calibrating such
models is a nuanced and well studied problem. For a guide to the literature on
the risks and returns of private equity investments, see
Kortweg~\cite{Korteweg2019}.
\paragraph{Continuous time.} There is a breadth of work on modeling portfolio
construction with illiquid assets. Many authors consider continuous time
stochastic processes. Dimmock et al.~study the endowment model, under which
university endowments hold high allocations in illiquid alternative assets, via
a continuous time dynamic choice model with deterministic-in-time discrete
liquidity shocks every $T$ periods~\cite{Dimmock2019}. They allow the investor to
increase the position in the illiquid asset instantaneously, not modeling the
delayed nature of capital calls. Ang et al.~also study a continuous time
problem, but model the timing of liquidity events of the illiquid asset as an
independent Poisson process~\cite{Ang2013}. Optimal solutions are assumed to have
almost surely non-negative liquid wealth, meaning that the investor must always
be able to cover the effects of illiquidity. Another important line of inquiry
studies the effects of illiquidity through the commitment risk of a fixed
alternative's commitment. Sorensen, Wang, and Yang~\cite{Sorensen2014} study this
problem by focusing on an investor who can modify their positions in stocks and
bonds, taking an investment in an illiquid asset as given and held to maturity.

\paragraph{Discrete time.} The discrete time case is also well studied.
Takahashi and Alexander first introduced what amounts to a deterministic linear system to
model an illiquid asset's calls, distributions, and asset
value~\cite{Takahashi2002}. This model posits that calls are a time-varying fraction of
uncalled commitments, and that distributions are a time-varying fraction of the illiquid
asset value, and returns are constant. Our model is similar, but differs in two important ways. 
First, our model is time-invariant. Second our model incorporates randomness in
these fractions as well as the returns. 
Giommetti and Sorensen use the Takahashi and Alexander model in a standard,
discrete-time, infinite-horizon, partial-equilibrium portfolio model to
determine optimal allocation to private equity~\cite{Giommetti2021}. Here the
calls and distributions are deterministic fractions of the uncalled commitments
and illiquid asset value, but the illiquid asset value grows with stochastic
returns. Again, there is an almost sure constraint which insists that the
investor's liquid wealth is never exhausted, which means that the investor
maintains a liquidity reserve of safe assets to cover calls.

\paragraph{Optimal allocation to illiquid assets.} Broadly, across the
literature we have reviewed, the reported optimal allocations to illiquid assets
are strikingly low compared to the de facto wants and need of institutional
investors who are increasingly allocating larger and larger shares of their
portfolios to illiquid alternatives. In their extensive survey of Illiquidity
and investment decisions, T\'edongap and Tafolong~\cite{Tedongap2018} report
that recommended illiquid allocations range from the low single digits to around
20\% on the upper end. This is strikingly lower than the target levels observed
in practice. For example, the National Association of College and University
Business Officers (NACUBO) provide data showing the allocation weights of
illiquid alternatives in University endowments reaching 52\% in 2010. Unlike
other analyses, our method does not require investors be able to cover calls
with probability one, and instead provides a tool for maintaining an optimized
target asset allocation under uncertain calls, distributions, returns, and
growth.

Hayes, Primbs, and Chiquoine propose a penalty cost approach to asset allocation
whereby an additional term is added to the traditional mean-variance
optimization (MVO) problem to compensate for the introduction of illiquidity
\cite{Hayes2015}. They solicit a user provided marginal cost curve which
captures the return premium needed for an illiquid asset to be preferred over a
theoretically equivalent liquid alternative. This leads to a formulation nearly
identical to the standard MVO problem, with a liquidity-adjusted expected return
(a function of the allocation). In their work the notion of liquidity is captured in a
scalar between $0$ and $1$.

\paragraph{Multi-period optimization.} Our policy is based on solving a
multi-period optimization problem. Dantzig and Infanger~\cite{Dantzig1993}
introduce a multi-stage stochastic linear programming approach to multi-period
portfolio optimization. Mulvey, Pauling, and Madey survey the advantages of
multi-period portfolio models, including the potential for variance reduction
and increased return, as well as the ability to analyze the probability of
achieving or missing goals \cite{Mulvey2003}. Boyd et al.~\cite{Boyd2017}
describe a general framework for multi-period convex optimization. This
framework focuses on planning a sequence of trades over a set of periods trades
given return forecasts, trading costs, and holding costs. Our framework also
solves a multi-period convex optimization problem, but we do not make an approximation
of the dynamics, which is more appropriate for the longer time horizons and thus
more significant growth observed in strategic asset allocation. 

\paragraph{Model predictive control.}
Our method falls
under the category of Model Predictive Control (MPC), which is both widely
studied in academia and used in industry. For a survey of MPC,
see for example the books \textit{Model Predictive Control}~\cite{mpcbook} or 
Garc\'ia et al.~\cite{Garcia1989}.
Herzog et al.~\cite{Herzog2006} use an MPC approach for 
multi-period portfolio optimization, but only consider normally distributed returns and 
standard liquid assets. They do include a factor model of returns, as well as a conditional
value at risk (CVaR) constraint which is different in interpretation but takes the same form
as our insolvency constraint. 
The closest work we have identified to our own is the thesis of Lee, who uses a
very similar multi-period optimization problem with linear illiquid dynamics
\cite{Lee2012}. We both use a quadratic risk, and use certainty equivalent
planning to solve an open loop control problem. Lee's problem is multi-period,
but the objective is a function of only the final period wealth, whereas in our
case we have stage costs, as well as constraints on the solvency of our
portfolio. Additionally, in our stochastic model we use random call and
distribution intensities.

\paragraph{Contributions.} The linear dynamics of the illiquid wealth motivate
model predictive control (MPC) as a solution method. To the best of our
knowledge, there do not exist multi-period optimization-based policies for
constructing portfolios with both liquid and illiquid alternative assets. We
believe our contributions include the following.
\begin{enumerate}
\item Incorporating random intensities with the classic linear model of the illiquid asset's calls and distributions.
\item Formulating a multi-period optimization problem to perform strategic asset allocation with liquid and illiquid assets.
\item Using homogeneous risk constraints to account for growth in the multi-period planning problem.
\item Using liquidity/insolvency constraints to ensure calls are covered with high probability.
\item Obtaining a performance bound for the problem by considering
a stylized liquid world where the illiquid asset is completely liquid.
\end{enumerate}

\clearpage
\section{Stochastic dynamic model for an illiquid asset}\label{stochastic-dynamic-model-illiquid}
In this section we describe our stochastic dynamic model of one illiquid asset.
Our model is closely related to the linear system proposed by Takahashi and Alexander~\cite{Takahashi2002}, with the addition of uncertainty in the capital calls
and distributions. We demonstrate a straightforward extension of our model which would include
the Takahashi model in \S\ref{s-extensions}. 
We consider a discrete-time setting, with period denoted by $t=1,2,3 \ldots$,
which could represent months, quarters, years, or any other period.
Our model involves the following quantities, all denominated in dollars.
\begin{itemize}
\item $I_t\geq0$ is the illiquid wealth (or position in or NAV 
of the illiquid asset) at period $t$.
\item $K_t\geq 0$ is the total uncalled commitments at period $t$.
\item $C_t\geq 0$ is the capital call at period $t$.
\item $D_t\geq 0$ is the distribution at period $t$.
\item $n_t\geq 0$ is the amount newly committed to the illiquid asset 
at period $t$.
\end{itemize}
The commitment $n_t$ is the only variable we can directly control or choose.
All the others are affected indirectly by $n_t$.

\paragraph{Dynamics.}  Here we describe how the variables evolve over time.
At period $t$,
\begin{itemize}
\item we make a new capital commitment $n_t$ (which we can choose) 
\item we receive capital call $C_t$ (which is not under our control)
\item we receive distribution $D_t$ (which is not under our control)
\end{itemize}
The uncalled commitment in period $t+1$ is
\[
K_{t+1}=K_t+n_t-C_t,
\]
and the illiquid wealth in period $t+1$ is
\[
I_{t+1}=I_tR_t+C_t-D_t,
\]
where $R_t\geq 0 $ is a random total return on the illiquid asset.

\paragraph{Calls and distributions.}
We model calls and distributions as random fractions of $K_t,P_t, $ and $n_t$.
We model calls as 
\[
C_t =\lambda_t^0n_t+\lambda_t^1K_t,
\]
where 
$\lambda_t^0$ $\in [0,1]$ is the random immediate commitment call intensity
and
$\lambda_t^1$ $\in [0,1]$ is the random existing commitment call intensity.
Similarly, we model distributions as
\[
D_t =I_tR_t\delta_t,
\]
where $\delta_t$ $\in [0,1]$ is the random distribution intensity.

We assume the random variables $(R_t,\lambda_t^0,\lambda_t^1,\delta_t) \in
\reals \times [0,1]^3$ are I.I.D.,
i.e., independent across time and identically distributed.
(But for fixed period $t$, 
the components $R_t$, $\lambda_t^0$, $\lambda_t^1$, and $\delta_t$ need not be
independent.)
We do not know these random variables when we choose the current commitment $n_t$.
Formally, we assume that $n_t \indep (R_t,\lambda_t^0,\lambda_t^1, \delta_t)$.
The current commitment can depend on anything known at the beginning of period $t$
(including for example past values of returns and intensities), but the current period
return and intensities are independent of the commitment.

\subsection{Stochastic linear system model}
The model above can be expressed as a linear dynamical system with random
dynamics and input matrices.  
With state $x_t=(I_t,K_t)\in\reals^2$ and the control or input
$u_t=n_t\in\reals$, the dynamics are given by
\[
x_{t+1}=A_tx_t+B_tu_t,
\]
where
\begin{equation}\label{e-illiquid-matrix-eqn}
A_t=\bmat{
R_t(1-\delta_t) &\lambda_t^1 \\
0 & 1-\lambda_t^1
},\qquad
B_t=\bmat{
\lambda_t^0\\
1-\lambda_t^0 
}.
\end{equation}
With output $y_t = (I_{t}, K_{t}, C_{t}, D_{t}) \in \reals^4$, we have
\[
y_{t}=F_tx_t+G_tu_t,
\]
where
\begin{equation}\label{e-output-eqn}
F_t=\bmat{
1 & 0\\
0 & 1\\
0 & \lambda^1_t\\
R_t(1-\delta_t) & 0
},\qquad
G_t=\bmat{
0\\
0\\
\lambda^0_t\\
0}.
\end{equation}

We assume the initial state is known.
We observe that $x_t \indep (F_t,G_t)$, since the former depends on
the initial state, $n_t$, and $(F_\tau,G_\tau)$ for $\tau <t$, and these
are all independent of $(F_t,G_t)$.

A careful reader might notice that these linear dynamics mean that the
commitments and distributions asymptotically approach zero but never terminate.
However, the fractions of calls and distributions relative to the initial
amounts are minuscule after several periods, and are negligible in the presence
of new commitments coming in each period. Additionally, Gupta and Van
Nieuwerburgh~\cite{Gupta2021} found in analyzing long-term private equity
behavior that often funds have activity even fifteen years after inception,
further justifying the persisting calls and distributions in the linear systems
model.

\subsection{Mean dynamics}\label{illiquid-mean-dynamics}
Let $\overline x_t = \Expect x_t$ denote the mean of the state, 
$\overline u_t = \Expect u_t$ denote the mean of the input or control, and
$\overline y_t = \Expect y_t$ denote the mean of the output.
We define the mean matrices
\[
\overline{A}=\Expect{A_t}, \qquad \overline{B}=\Expect{B_t}, \qquad
\overline{F}=\Expect{F_t}, \qquad \overline{G}=\Expect{G_t}
\]
(which do not depend on $t$).
We then have
\begin{equation}\label{e-mean-dynamics}
\overline{x}_{t+1} = \overline{A}\overline{x}_t+\overline{B}\overline{u}_t, \qquad 
\overline{y}_{t} = \overline{F}\overline{x}_t+\overline{G}\overline{u}_t,
\end{equation}
which states that the
mean state and output is described by the same linear dynamical system, with
the random matrices replaced with their expectations. 
The mean dynamics is a time-invariant deterministic linear dynamical system.

\subsection{Impulse and step responses}

Our linear system model implies that
the mapping from the sequence of commitments or inputs $u_1, u_2, \ldots$ to the 
resulting outputs $y_1, y_2, \ldots$ is linear but random.
We review three standard concepts from linear dynamical systems.

\paragraph{Commitment impulse response.} 
We can consider the response of uncalled commitments, calls, illiquid wealth,
and distributions to committing $n_1=1$ at period $t=1$ and $n_t=0$ for all 
$t>1$, with zero initial state.
This is referred to as the impulse response of the system. The impulse
response is the stochastic process
\[
y_t=\bmat{
I_t\\
K_t\\
C_t\\
D_t
}= F_tA_{t-1}\cdots A_2 B_1, \quad t=1,2, \ldots .
\]
From the mean dynamics \eqref{e-mean-dynamics}, we know that the
mean impulse response is given by
\[
\overline{y}_t = \overline{F}~{\overline A}^{t-1}\overline B , \quad t=1,2, \ldots .
\]

\paragraph{Commitment step response.} We can also consider the effect of 
committing $n_1=n_2 = \cdots = 1$,
which is referred to as the step response of the
system. The step response shows how a simple policy of constant commitment
impacts our exposure over time in the illiquid asset, calls, distributions, and our level
of uncalled commitment. 
The step response is the stochastic process 
\[
y_t = F_t\paren{\paren{\sum_{i=0}^{t-3}A_{t-1}\cdots A_{2+i}B_{i+1}}+B_{t-1}}+G_t,\quad t=1,2, \ldots.
\]
From the mean dynamics \eqref{e-mean-dynamics}, we know that the
mean step response is given by
\begin{equation}\label{e-mean-step-response}
\overline y_t = \overline F\paren{\sum_{i=0}^{t-2}\overline A^i}\overline B+\overline G,\quad 
t=1,2, \ldots .
\end{equation}

\paragraph{Steady state response.}
We define the unit steady-state mean response $y^\text{ss}$ as
$\lim_{t \to \infty} \overline  y_t^\text{step}$.
Assuming the spectral radius of $\overline A$ is less than one, 
we take the limit of \eqref{e-mean-step-response} as $t\to\infty$ to obtain 
\[
y^\text{ss} = \overline F(I-\overline A)^{-1} \overline B + \overline G.
\]
We refer to the entries of 
\begin{equation}\label{e-steady-state-gains}
y^\text{ss} = (\alpha_K,\alpha_C,\alpha_I, \alpha_D)
\end{equation}
as the steady state gains from commitment to illiquid wealth, uncalled commitment,
capital calls, and distribution.
These numbers have a simple interpretation.  
For example, $\alpha_I$ tells us what the asymptotic mean illiquid wealth is, 
if we repeatedly make a commitment of \$1.
It can be shown that $\alpha_C=1$, \ie, if we constantly commit \$1, then 
asymptotically, and in mean, the capital calls will also be \$1.

\subsection{A particular return and intensity distribution}
We suggest the following parametric joint distribution for
$(\lambda_t^1,\lambda_t^0,\delta_t,R_t)$. They are generated from a random
3-vector
\begin{equation}\label{e-illiquid-latent}
z_t\sim \mathcal{N}(\mu,\Sigma)\in\reals^3.
\end{equation}
From these we obtain 
\begin{equation}\label{e-illiquid-intensity-model}
\lambda_t^1= \frac{1}{1+\exp(z_t)_1},\qquad
\lambda_t^0 = \frac{1}{2}\lambda_t^1,\qquad
\delta_t = \frac{1}{1+\exp {(z_t)}_2},\qquad
R_t= \exp {(z_t)}_3.
\end{equation}
With this model, the return is log-normally distributed while the call and distribution 
intensities are logit-normally distributed.  Dependency among the return and 
the intensities are modeled by the off-diagonal entries of $\Sigma$.

\subsection{Example}\label{example_illiquid_stochastic}
Here we describe a particular instance of the distribution described above,
that we will use in various numerical examples in the sequel.

\paragraph{Example return and intensity distribution.} In this example we use
the following parameters for the distribution of
$(\lambda_t^1,\lambda_t^0,\delta_t,R_t)$ specified in \eqref{e-illiquid-latent}:
\begin{equation}\label{e-intensity-dist}
\mu = \bmat{-0.700\\-0.423\\ 0.158},\qquad
\Sigma = \bmat{
0.068   & 0.072 & 0.006\\
0.073   & 0.271 & 0.043\\
0.006   & 0.043 & 0.079
}.
\end{equation}

This example is based on yearly periods.
The mean return of the illiquid asset is derived from the BlackRock Capital Market
Assumptions as of July 2021, which reports one private equity asset,
Buyout, with a mean annual return of 15.8\% \cite{cma}. 
The call and distribution mean intensities are calibrated from private equity data
for the eFront Buyout fund. 
The mean values of the intensities  are 
we report the empirical means 
\[
\overline \lambda_t^1 = .26,\qquad \overline \lambda_t^0 = .128,\qquad \overline \delta_t=.33.
\]
(These are found by Monte Carlo simulation, since
the mean of a logit-normal distribution doesn't have an analytical expression.)
The covariance matrix is calibrated from the same data. 

\begin{figure}
\centering
\includegraphics[scale=.45]{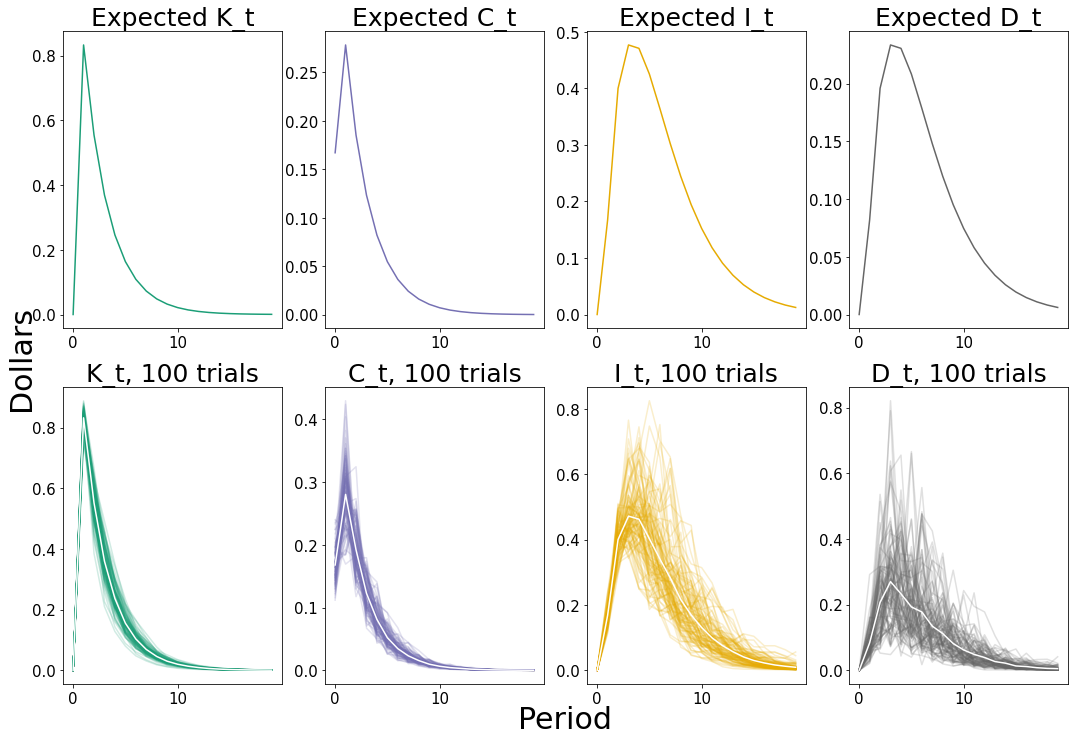}
\caption{Commitment impulse response. 
The top plot shows the mean values,
and the bottom plot shows 100 realizations of the stochastic model.}
\label{impulse_response}
\end{figure}

\paragraph{Commitment impulse response.} The impulse response from commitment to
uncalled commitment, calls, illiquid wealth, and distributions is shown in
figure~\ref{impulse_response}. The top row shows the mean response, and the
bottom shows 100 realizations, with the empirical mean shown as the white curve.

We see that the uncalled commitments peak in the next period at a level of about
0.8. The calls peak at the next period at .28. We can see that our initial
commitment translates into an illiquid holding which, in expectation, peaks four
periods later with a value of about .47. Similarly, the resulting distributions
peak with the illiquid wealth four periods later, with a level of .24.

\begin{figure}
\centering
\includegraphics[scale=.45]{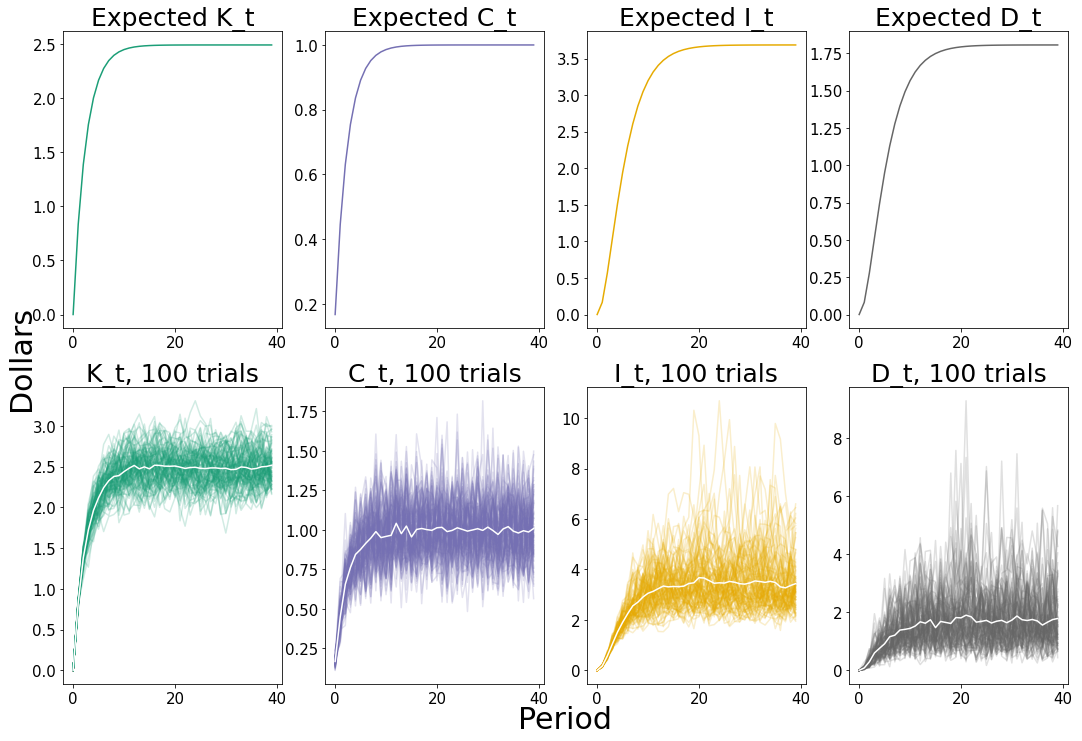}
\caption{Commitment step response.
The top plot shows the mean values,
and the bottom plot shows 100 realizations of the stochastic model.}
\label{step_response}
\end{figure}

\paragraph{Commitment step response.} In figure~\ref{step_response}, we see the
step response to constant unit commitment of uncalled commitments, calls,
illiquid wealth, and distributions. The top row shows the mean response, and the
bottom shows 100 realizations.
The mean step responses converge in around 8 periods to values near their 
asymptotic values.

\paragraph{Asymptotic expected response to constant commitment.}
The steady-state gains are
\[
\alpha_K= 2.491, \qquad
\alpha_C= 1.000, \qquad
\alpha_I= 3.685, \qquad
\alpha_D= 1.804.
\]
For example, repeatedly committing \$1 leads to an asymptotic mean illiquid wealth
of \$3.685.

\subsection{Comparison with the Takahashi and Alexander model}
Our stochastic model of an illiquid asset is closely related to that of Takahashi
and Alexander \cite{Takahashi2002}, but it differs in to key ways. The most important
difference is that our model is Markovian; the calls, distributions, and 
returns at time $t$ are conditionally independent of the all previous quantities, 
given the state at time $t$. In comparison, Takahashi and Alexander's model specifies
time varying call and distribution intensity parameters. 
These time varying intensities mean that the final intensities can be set to 1, 
meaning calls and distributions can have deterministic end times, and the exposure 
will not geometrically decline. In \S\ref{s-extensions} we describe how to modify our model
to depend on arbitrarily many previous time periods. 
This means we can exactly capture the original Takahashi and Alexander model with this extension of our model.
We emphasize that this generalization remains fully tractable from the portfolio optimization standpoint
described in this paper. The second difference between our model and that of Takahashi and Alexander is
that ours is a stochastic model, with random intensities, whereas theirs is deterministic. 
\clearpage

\section{Commitment optimization}\label{s-commitment-opt}

Because the dynamics is linear, we can use
convex optimization to choose a sequence of commitments to meet various goals
in expectation.
We consider a simple example here to illustrate.

We consider the task of starting with no illiquid exposure or uncalled commitments,
\ie, $I_1=0$, $K_1=0$, and choosing a sequence of commitments, $n_t$, $t=1,\ldots,T$.
Our goal is to reach and maintain an illiquid wealth of $I^\text{tar}$,
so we use as our primary objective the mean-square tracking error,
\[
\Expect  \frac{1}{T+1} \sum_{t=1}^{T+1} (I_t - I^\text{tar})^2.
\]
In addition, we want a smooth sequence of commitments, so we add a secondary objective
term which is the mean square difference in commitments,
\[
\Expect \frac{1}{T-1} \sum_{t=2}^T (n_t-n_{t-1})^2,
\]
We also impose a maximum allowed per-period commitment, \ie, $n_t \leq n^\text{lim}$.
(Of course we can add other constraints and objective terms; this example is meant only to
illustrate the idea.)

We illustrate two methods.  The first is simple open-loop planning, 
in which assume that state follows its mean trajectory, and we determine a fixed
sequence of commitments, and then simply execute this plan.
The second method is a closed-loop method, which adapts the commitments based
on previously realized returns, capital calls, and distributions.
This method is called model predictive control (MPC).
We evaluate both policies under the mean dynamics and the stochastic dynamics.

\subsection{Open loop commitment control} 
\paragraph{Planning.} We will find a plan of commitments based on the mean dynamics.
This leads to the convex optimization problem
\[
\begin{array}{ll}
  \mbox{minimize} & \frac{1}{T+1}\sum_{t=1}^{T+1}(\hat I_t-I^\text{targ})^2
  +\gamma^\text{smooth}\frac{1}{T-1}\sum_{t=2}^T(\hat n_t-\hat n_{t-1})^2\\
  \mbox{subject to} & \hat x_{t+1}=\overline{A}\hat x_t
  +\overline{B}\hat n_t,\qquad t=1,\ldots,T\\
  & \hat x_1 = 0\\
  & 0\leq \hat n_t\leq n^\text{lim},\qquad t=1,\ldots,T,
\end{array}
\]
where $\gamma^\text{smooth}>0$ is a hyperparameter that determines the weight of
the smoothing penalty, and $\overline A$ and $\overline B$
are as defined in \eqref{e-output-eqn}. The variables in this problem are $\hat n_1,\ldots,\hat n_T$ 
and $\hat x_1,\ldots,\hat x_{T+1}$, with $\hat I_t=(x_t)_1$ for $t=1,\ldots,T+1$. 
We use the notation $\hat I_t$, $\hat n_t$, $\hat x_t$ to emphasize that these are quantities in our
plan, and not the realized values.
This is a simple convex optimization problem, a quadratic program (QP), and readily solved 
\cite{cvxbook}.

\begin{figure}\label{commit_plan}
\begin{center}
  \includegraphics[width=0.8\textwidth]{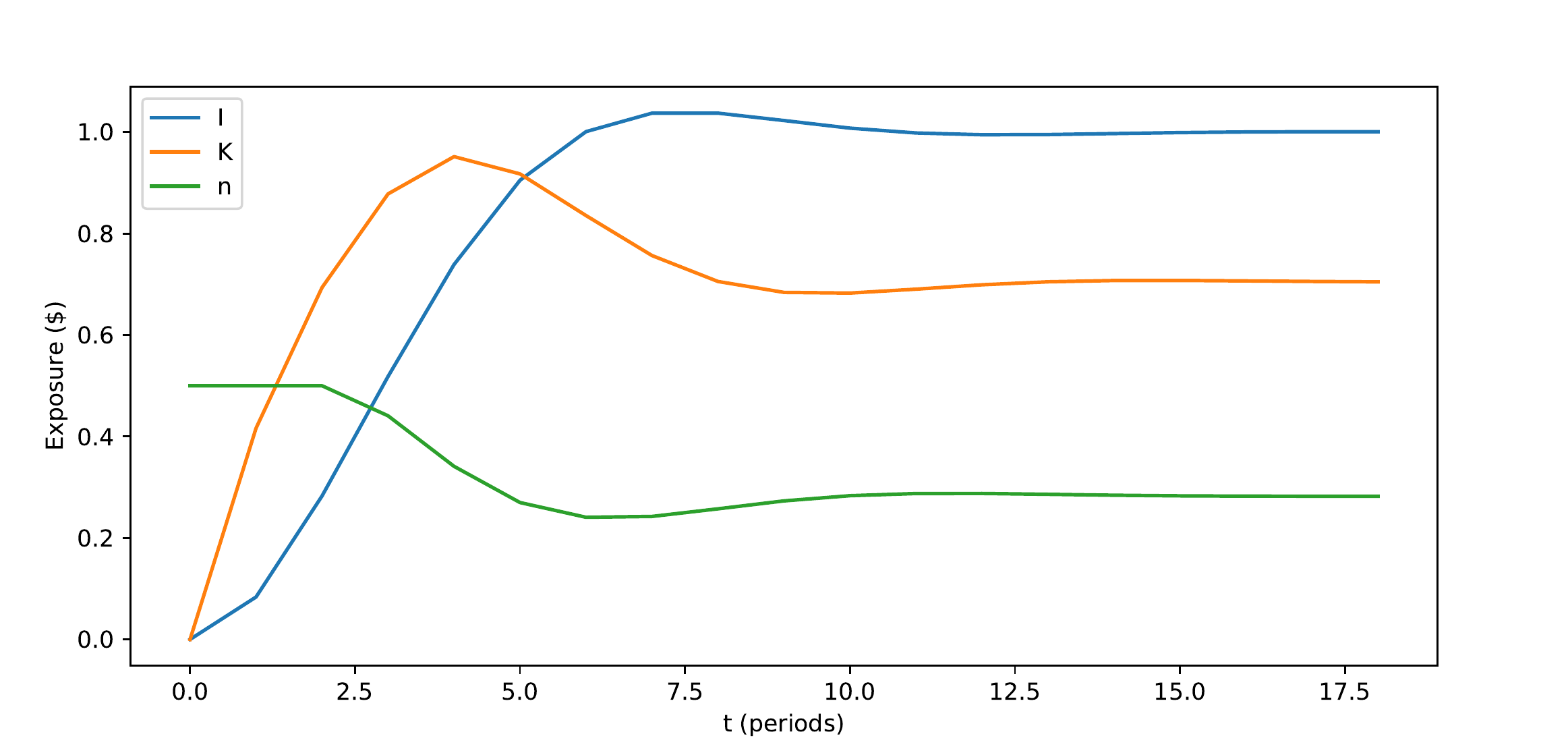}
\end{center}
\caption{Deterministic commitment planning.}
\end{figure}

\paragraph{Example.} Solving this problem with our running example defined in \eqref{e-intensity-dist} for
\begin{equation}\label{e-commitment-planning-params}
T=20,\qquad \gamma^\text{smooth}=1,\qquad I^\text{targ}=1,\qquad n^\text{lim}=.5,
\end{equation}
we find the optimal planned sequence of commitments and
corresponding uncalled commitment and exposure shown in
figure~\ref{commit_plan}.

The mean-squared tracking error attained by our plan is 0.133. We can also calculate the tracking error from $t=5$ onwards to account for the large contribution
to tracking error of the first four periods.
Thus, a perhaps more meaningful metric is the delayed root-mean-square (RMS) tracking error, 
\[
\paren{\frac{1}{T-4}\sum_{t=5}^{T}(I_t-I^\text{targ})}^{1/2}.
\]
(This is on the same scale as $I^\text{targ}$ and is readily compared to it.)
The plan attains a delayed RMS tracking error of 0.071.

The optimal commitment sequence makes sense.  It hits the limit for the first two
periods, in order to quickly bring up the illiquid wealth; then it backs off to a lower 
level by around period $6$, and finally converges to
an asymyptotic value near $I^\text{targ}/\alpha_I = 0.27$, which is the constant commitment
value that would asymptotical give mean illiquid value $I^\text{targ}$.

\begin{figure}
\begin{center}
\includegraphics[width=0.8\textwidth]{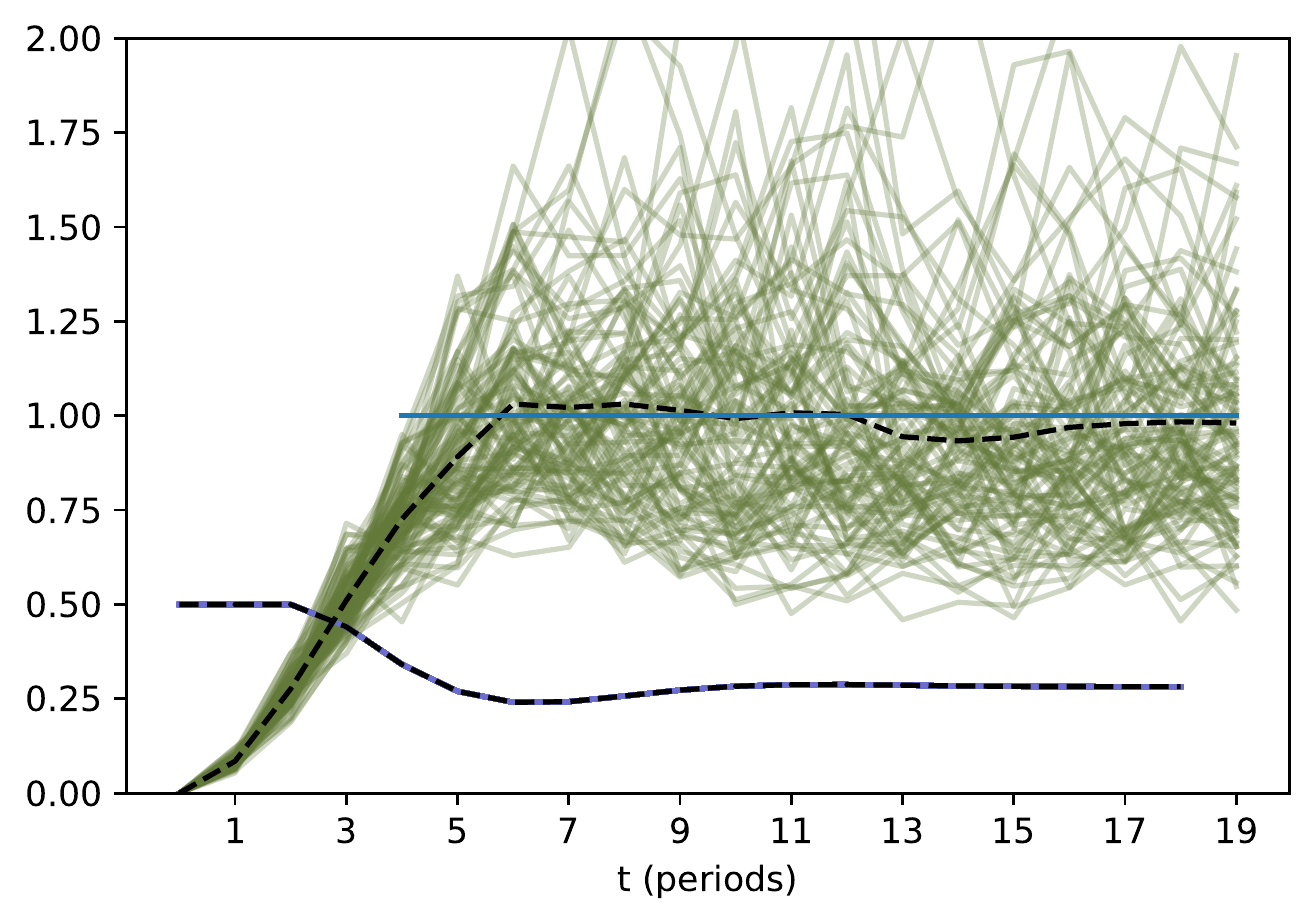}
\end{center}
\caption{Open loop commitment control, 100 realizations.}
\label{commit_open}
\end{figure}

\paragraph{Results.}
These results above are with the mean dynamics.
We can also execute this sequence of planned commitments under random calls and 
distributions as
specified by our stochastic model in \S\ref{example_illiquid_stochastic}.
The results for 100 simulated realizations is shown in figure \ref{commit_open}.
The mean-squared tracking error, averaged across the realizations, is 0.199.
The delayed root-mean-square tracking error, averaged across the realizations
is 0.274.

\subsection{Closed loop commitment control via MPC} 
We now perform model predictive control, which is a \emph{closed loop} method, 
meaning $n_t$ can depend on $x_t$, \ie, we can adapt our commitments to the
current values of uncalled commitments and illiquid wealth.

\begin{figure}
\begin{center}
\includegraphics[width=0.8\textwidth]{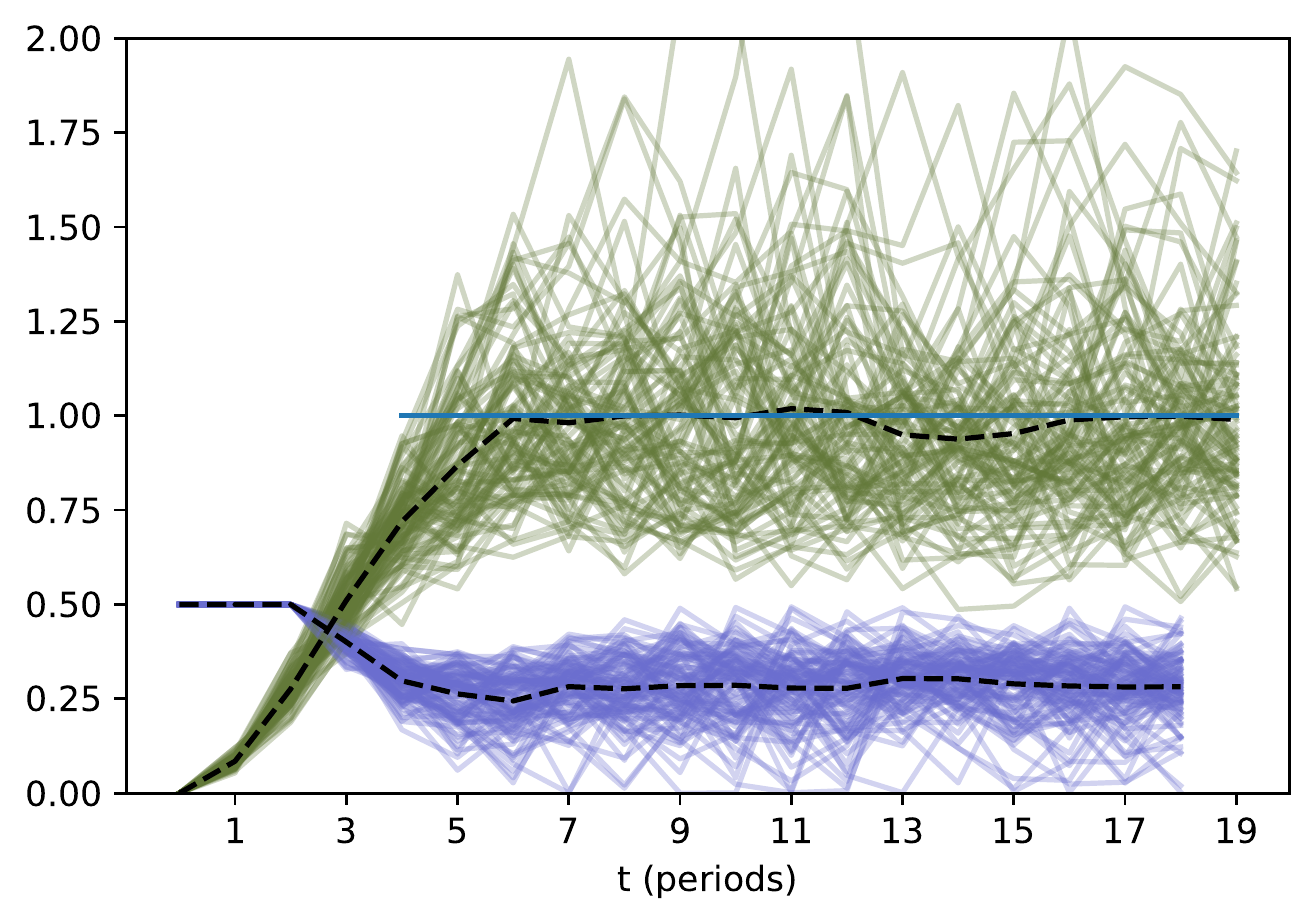}
\end{center}
\caption{MPC commitment control, 100 realizations.}
\label{commit_closed}
\end{figure}

\paragraph{Planning.} At every time $t=1,\ldots,T$, 
we plan commitments over the next $H$ periods $t,t+1,\ldots,T+H$,
where $H$ is a planning horizon.
We use $\hat x_{\tau|t},\hat n_{\tau|t}$
to indicate that these are the quantities in the plan at time $\tau$, from the plan made at time $t$. 
These planned quantities are found by solving the optimization problem
\[
\begin{array}{ll}
  \mbox{minimize} & \frac{1}{H+1}\sum_{\tau=t}^{t+H+1}(\hat I_{\tau|t}-I^\text{targ})^2+
  \gamma^\text{smooth}\frac{1}{H-1}\sum_{t=2}^T(\hat n_{\tau|t}-\hat n_{\tau-1|t})\\
  \text{subject to} & \hat x_{t|t}=x_t\\
  & \hat x_{\tau+1|t}=\overline{A}\hat x_{\tau|t}+\overline{B}\hat n_{\tau|t},\quad\tau=t,\ldots,t+H\\
  & 0\leq \hat n_{\tau|t}\leq n^\text{lim},\quad \tau=t,\ldots,t+H,
\end{array}
\]
with variables $\hat x_{t|t},\ldots,\hat x_{t+H+1|t}$ and $ \hat n_{t|t},\ldots,\hat n_{t+H|t}$,
$\hat I_{\tau|t} = (\hat x_{t|t})_1$
Note that when we plan at time $t$, we include the constraint $\hat x_{t|t} = x_t$; 
this closes the feedback loop by planning based on the current realized state.
\paragraph{Policy.} Our policy is simple to explain: at time $t$, 
after planning as described above, we execute control
\[
n_t=\hat n_{t|t}.
\]
Note that the planned quantities $\hat I_{\tau|t}$, $\hat x_{\tau|t}$, $\tau=t+1,\ldots,t+H+1$, 
and $\hat n_{\tau|t}$, $\tau=t+1,\ldots,t+H$, are never executed by the MPC policy. 
They are only part of the plan.
\paragraph{Results.} We again execute our policy under random calls and distributions as
specified by our stochastic model in \S\ref{example_illiquid_stochastic}.
The results for 100 simulated trajectories a shown in figure \ref{commit_closed}.
The average mean-squared error is 0.182.
The average delayed root-mean-square tracking error is 0.244, 
an 11\% reduction from the open loop policy.

\clearpage
\section{Joint liquid and illiquid model}\label{joint_illiquid_liquid_model}
We now describe a model for an investment universe consisting of multiple illiquid alternative
and liquid assets. 
First, we extend to a universe of $n^\text{ill}$ illiquid assets. 
\paragraph{Multiple illiquids.} 
We extend the same quantities as in \S\ref{stochastic-dynamic-model-illiquid} from scalars
to vectors of dimension $n^\text{ill}$.
\[
K_t,I_t,C_t,D_t \in \reals^{n^\text{ill}},
\qquad n_t\in\reals^{n^\text{ill}},
\qquad R_t^\text{ill},\lambda_t^1,\lambda_t^0,\delta_t \in\reals^{n^\text{ill}}.
\]
We have the exact same dynamics as before, duplicated for each illiquid asset.
Each has its own states for exposure and uncalled commitment, and its own control 
for its new commitments. 
The illiquid calls, distributions, and returns are now part of a joint distribution. 
The illiquid dynamics extend in vectorized form to 
\begin{equation*}
K_{t+1} = K_t+n_t-C_t,\qquad I_{t+1} = \diag(R_t)I_t+C_t-D_t,
\end{equation*}
with
\begin{equation*}
C_t = \diag(\lambda_t^0)n_t+\diag(\lambda_t^1)K_t, 
\qquad D_t = \diag(R_t)\diag(\delta_t)I_t.
\end{equation*}
We emphasize that while the return, call, and distribution dynamics here are separable
across the illiquid assets, the underling random variables
$((R_t)_j,(\lambda_t^0)_j,(\lambda_t^1)_j,(\delta_t)_j)$
can be modeled jointly.
We continue with our assumption that these random variables are independent across time.

\paragraph{Multiple liquids.}
There are now a set of $n^\text{liq}$ liquid assets available to us. 
The liquid assets are simple: we can buy and sell them at will at each period;
they suffer none of the complex dynamics of the illiquid assets.
We add one new state, $L_t$, the (total) liquid wealth at period $t$.
In addition to new commitments for each illiquid asset, at each time $t$ we now control how we allocate our liquid wealth
each period, as well as how much outside cash to inject into our liquid wealth. 
Thus we have the additional quantities, which we can control:
\begin{itemize}
\item $h_t\geq 0\ (\in \reals^{n^\text{liq}})$ is the allocation in dollars invested in liquid assets at period $t$
\item $s_t\geq 0\ (\in \reals)$ is the outside cash injected at period $t$
\end{itemize}
At the beginning of period $t$, we invest (or allocate) our liquid wealth in liquid assets.
This corresponds to the constraint $L_t=\mathbf{1}^Th_t$. 
We receive multiplicative liquid returns $(R_t^\text{liq})_j\in \reals$ on liquid asset $j$,
yielding total return $h_t^T R_t^\text{liq}$.
We then pay out capital calls from and receive distributions to our liquid wealth, for all illiquid assets. 
This corresponds to a net increase in liquid wealth given by $-\ones^TC_t+\ones^TD_t$. 
Lastly, if at this stage our liquid wealth is negative, we are forced to add outside cash $s_t$
to at least bring our liquid wealth to zero.
Compactly, the liquid dynamics are
\[
L_{t+1}=h_t^T R_t^\text{liq}-\ones^TC_t+\ones^TD_t+s_t,
\]
with constraints
\begin{equation*}
  h_t,n_t,s_t\geq 0, \qquad L_t\geq 0, \qquad L_t=\mathbf{1}^Th_t.
  \end{equation*}

\subsection{Stochastic linear system model}\label{stochastic_linear_system_joint}

We can again represent these dynamics as a stochastic linear system. 
Let $x_t = (L_t,I_t,K_t)\in \reals^{1+2n^\text{ill}}$ be the state vector. 
The control is $u_t = (h_t,n_t,s_t)\in \reals^{1+n^\text{liq}+n^\text{ill}}$.
Extending the $A$ and $B$ matrices from \S\ref{e-illiquid-matrix-eqn}, define
\begin{equation}\label{e-joint-matrix-dynamics}
A_t = \bmat{
0   & (\delta_t\circ R_t^\text{ill})^T      & {-\lambda^1_t}^T\\
0   & \diag((\ones-\delta_t)R_t^\text{ill}) & \diag(\lambda_t^1)\\
0   & \diag(\zeroes)                        & \diag(\ones-\lambda_t^1)
},\qquad
B_t = \bmat{
{R^\text{liq}_t}^T  & {-\lambda_t^0}^T          &  1\\ 
\zeroes^T           & \diag(\lambda_t^0)        &  0\\
\zeroes^T           & \diag(\ones-\lambda_t^0)  &  0}.
\end{equation}
Then the random linear dynamics with multiple illiquids and liquids are
\[
x_{t+1}=
A_t x_t+ B_tu_t,
\]
with constraints
\begin{equation*}
h_t,n_t,s_t\geq 0, \qquad L_t\geq 0, \qquad L_t=\mathbf{1}^Th_t.
\end{equation*}
The presence of the outside cash control $s_t$ implies that a feasible control exists
for any feasible value of the states, since $s_t$ prevents the liquid wealth from ever being negative.

As in \S\ref{illiquid-mean-dynamics}, we let $\overline x_t = \Expect x_t$ denote the mean of the state, 
$\overline u_t = \Expect u_t$ denote the mean of the input or control, define the mean system matrices as
\begin{equation*}
\overline{A}=\Expect{A_t},\qquad \overline{B}=\Expect{B_t},
\end{equation*}
and recover the same mean dynamics
\begin{equation}\label{e-mean-liquid-dynamics}
\overline{x}_{t+1} = \overline{A}\overline{x}_t+\overline{B}\overline{u}_t.
\end{equation}
\subsection{Return and intensity distribution}
We extend the previous generative model specified in \eqref{e-illiquid-latent}
and \eqref{e-illiquid-intensity-model} to include liquid returns, 
\begin{equation}\label{e-join-latent}
z_t=\bmat{
z_t^\text{int}\\
z_t^\text{ret}
}
\sim \mathcal{N}(\mu,\Sigma)\in\reals^{3n^\text{ill}+n^\text{liq}},\qquad
\mu=\bmat{\mu^\text{int}\\
\mu^\text{ret}\\},\qquad 
\Sigma=\bmat{
\Sigma^\text{int} & \Sigma_{12}\\
\Sigma_{21} & \Sigma^\text{ret}
}.
\end{equation}
From these we obtain delayed and immediate call intensities
\[ 
\lambda_t^1= \frac{1}{1+\exp(z_t)_{1:n^\text{ill}}},\qquad
\lambda_t^0 = \frac{1}{2}\lambda_t^1,\]
distribution intensities
\[
\delta_t = \frac{1}{1+\exp+(z_t)_{n^\text{ill}+1:2n^\text{ill}}},
\]
and returns
\[
\bmat{
R_t^\text{ill}\\
R_t^\text{liq}
}
= \bmat{
\exp {(z_t)}_{2n^\text{ill}+1:3n^{n^\text{ill}}}\\
\exp {(z_t)}_{3n^\text{ill}:}}.
\]

\clearpage
\section{Strategic asset allocation under the relaxed liquid model}\label{relaxed-SAA}

In this section we introduce a highly simplified model, where all of the challenges of
illiquid alternative assets are swept under the rug.
This model is definitely not realistic, but we can use it to develop an unattainable benchmark
for performance that can be obtained with the more accurate model.

\subsection{Relaxed liquid model}\label{Relax}

As a thought experiment, we
imagine the illiquid assets are completely liquid: we have arbitrary control of
illiquid asset positions (immediate increase or decrease).  This is a relaxation
of the actual problem setting, where we must face stochastic and only indirectly
controllable calls and distributions.  The idea of a relaxed liquid model is not new;
for example, Giommetti et al.~\cite{Giommetti2021}
consider the target allocations resulting from treating illiquid assets as fully
liquid for comparison, but do not evaluate stochastic control policies trying to
achieve these allocations in an illiquid world.
The relaxed liquid model is also implicitly behind various Captial Market Assumptions,
where return ranges, and correlations, are given for both liquid and illiquid assets.

The relaxed liquid model is very simple.  There is only one state, the total wealth $W_t$.
The quantities we have control over are the allocations to liquid and illiquid assets,
denoted $h_t^\text{liq} \in \reals^{n^\text{liq}}$ and
$h_t^\text{ill} \in \reals^{n^\text{ill}}$.
The wealth evolves according to the dynamics
\[
W_{t+1} = u_t^Tr_t,\qquad
u_t^T\ones=1,\qquad
u=\bmat{
h_t^\text{liq}\\
h_t^\text{ill}
},
\]
where $r_t=z_t^\text{ret}$ is defined in \eqref{e-join-latent}.

We use the standard trick of working with the weights of the allocations
in each period, denoted $w_t$, instead of $u_t$.  This is defined as 
$w_t = u_t/W_t$, so $\ones^T w_t=1$.  We recover the dollar allocations as
$u_t = W_tw_t$.

\subsection{Markowitz allocation and policy}\label{markowitz}

A standard way to choose a portfolio allocation is to
solve the one period risk-constrained Markowitz problem,
\begin{equation}
\begin{array}{ll}\label{e-markowitz-problem}
\text{maximize} & \mu^Tw\\
\text{subject to} &\ones^Tw =1,\quad w\geq 0\\
& \|\Sigma^{1/2}w\|_2 \leq \sigma,\\
\end{array}
\end{equation}

where $\sigma$ is the maximum tolerable return standard deviation, and $\mu$ and
$\Sigma$ are the expected return and return covariance, respectively. 
We denote the optimal allocation as $w^\star$.
The natural policy associated with solving the Markowitz
problem simply rebalances to $w^\star$: 
it sets $u_t=W_tw^\star$ for each period $t$.
This simple rebalancing is of course not possible under the accurate model that includes
the challenges of alternative assets, but it is under the relaxed liquid model.

\subsection{Example}\label{relaxed_example}

\paragraph{Liquid performance.} Under the assumptions of \S\ref{Relax}, we solve
the one period Markowitz problem with $\mu_\text{ret},\Sigma_\text{ret}$ the
mean and covariance of the joint distribution of liquid, illiquid asset returns.
Using this relaxed Markowitz target, we simulate the fantasy performance
achieved by being able to perfectly rebalance both liquids and illiquids to the
Markowitz target each period, for multiple periods, using the policy described
earlier in \ref{markowitz}. 

For the parameters defined in
\eqref{e-join-latent}, we use the specific values of
{
\begin{eqnarray}\label{e-example-mean-return}
\mu_\text{ret} &=& \pmat{
0.158&
0.000&
0.072&
0.023&
0.036&
0.046},\\
\sigma_\text{ret} &=& \pmat{
0.281&
0.000&
0.206&
0.046&
0.047&
0.162},\\
C_\text{ret} &=& \bmat{
1.000  & 0.000  & 0.422  & -0.298 & -0.002 & 0.261  \\
0.000  & 1.000  & 0.000  & 0.000  & 0.000  & 0.000 \\
0.422  & 0.000  & 1.000  & -0.843 & 0.197  & 0.800  \\
-0.298 & 0.000  & -0.843 & 1.000  & -0.018 & -0.739 \\
-0.002 & 0.000  & 0.197  & -0.018 & 1.000  & 0.628  \\
0.261  & 0.000  & 0.800  & -0.739 & 0.628  & 1.000  
},\label{e-example-cov-return}
\end{eqnarray}
with $\Sigma_\text{ret} = \diag(\sigma_\text{ret})C\diag(\sigma_\text{ret})$.
The liquid return mean and covariance matrix are gathered from the BlackRock
Capital Market Assumptions for equities as of July 2021 \cite{cma}.
The corresponding expected returns are
\[
\overline R_t^\text{liq}= \pmat{1.21834088 & 1.0976099 & 1.02436202 & 1.0377968 & 1.0609715 & 1},
\]
where the last asset is cash. 

\paragraph{Risk-return trade-off.} By solving the Markowitz problem with these
parameters across a range of values for the risk tolerance $\sigma$ (which give
rise to corresponding Markowitz targets), we can create a risk-return trade-off
plot, shown in figure~\ref{baseline_results}.
We should consider this trade-off curve as an unattainable performance benchmark, that we 
can only strive to attain when the challenges of illiquid alternatives are present.
\begin{figure}
\begin{center}
\includegraphics[scale=.3]{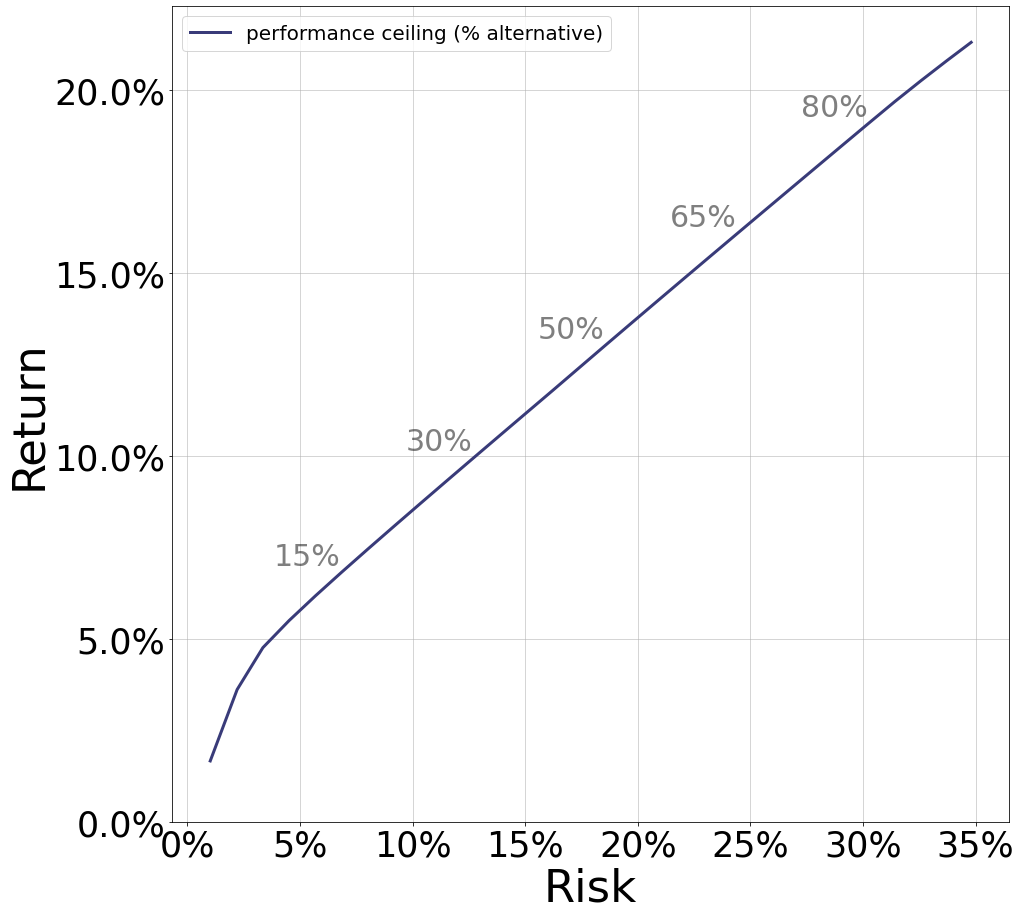}
\caption{Relaxed liquid risk-return trade-off, obtained by pretending the illiquid assets
are fully liquid. This gives an unattainable performance benchmark for problem when
the challenges of illiquid alternatives are present.}
\label{baseline_results}
\end{center}
\end{figure}


\clearpage
\section{Strategic asset allocation with full illiquid dynamics} 

In \S\ref{relaxed-SAA}, we describe an approach to strategic asset allocation
for portfolios including an imagined class of illiquid alternatives which are
rendered completely liquid. In this section, we provide methods to perform 
strategic asset allocation with mixed liquid and illiquid alternative portfolios
where we can only augment our illiquid position by making new commitments, 
and the effect of this action is random and delayed. First, we describe a method
which over time establishes and then maintains a given target allocation under growth.
Then, we describe a more sophisticated MPC method which jointly selects a target
allocation based on a user's risk tolerance, establishes the target, and maintains the
target in growth.

\subsection{Steady-state commitment policy}\label{heuristic}

We first describe a simple policy, which seeks to track a target allocation $\theta^\text{targ}$.
It allocates liquid assets proportionate to its desired liquid allocation, and
makes new commitments of a target level of illiquid wealth scaled by the
asymptotic expected private response to constant commitment. The input is a
target allocation $\theta^\text{targ}$, current liquid wealth L and illiquid
wealth I. First, the policy checks if $L$ is negative. If it is, it returns
control \[ u=(h,n,s),\qquad h=0,\quad n=0,\quad s=|L|.  \] Otherwise, if the
liquid wealth is positive, the policy proceeds as follows.  First, the policy
rebalances the liquid holdings proportionately to $\theta^\text{targ}$,
\[
h=L\frac{\theta^\text{liq}}{1^T\theta^\text{liq}},
\]
where $\theta^\text{liq},\theta^\text{ill}$ are the liquid and illiquid blocks of
the allocation vector 
$\theta^\text{targ}=\bmat{\theta^\text{liq}\\\theta^\text{ill}}$, respectively.
Then,  with $\alpha_I$ as the 1 dollar private commitment step response defined
in \eqref{e-steady-state-gains}, and $I^\text{targ}$ as the target illiquid level, 
$p^\text{targ}=\theta^\text{ill}(L+I)$, the policy commits 
\[
n_i=\frac{I^\text{targ}_i}{\alpha_{I_i}}
\]
and returns control $u=(h,n,0)$.

\subsection{Model predictive control policy}\label{MPC} We now describe a more
sophisticated policy which plans ahead based on a model of the future, seeking
to maximize wealth subject to various risk constraints.  For a sequence of
prospective actions, the policy forecasts future state variables using the mean
dynamics described in \eqref{e-mean-liquid-dynamics}. 
The policy then chooses a sequence of actions by optimizing an
objective which depends on the planned actions and forecast states. Finally,
the policy executes solely the first step of the planned sequence. The impact
of that action is observed, and the resulting state is observed, and then this
cycle repeats. 

The policy selects a planned sequence of actions by trying to maximize the ultimate
total liquid and illiquid wealth.
However, it is also constrained by a user's risk tolerance, which caps the 
allowable per period return volatility. 
Additionally, because capital calls are stochastic in nature, the policy seeks to guarantee
that with high probability, all capital calls can be funded from the liquid wealth.

\paragraph{Modified Markowitz constraint.} Motivated by the standard one period
risk-constrained Markowitz problem \eqref{e-markowitz-problem}, we would like to include a risk constraint
in our planning problem.  However, the Markowitz problem has variables in
weight space rather than wealth space.  Other multi-period optimization
problems based on the Markowitz problem, such as in \cite{Boyd2017}, assume a
timescale over which the wealth does not grow significantly over the planning horizon.
In our case, since potential application contexts include endowments and
insurers, we must handle substantial growth over the investment horizon. 
Thus, we consider an analogous risk
constraint in wealth space rather than weight space,
\[
\frac{y^T\Sigma y}{(\ones^Ty)^2}\leq \sigma^2 \iff \|\Sigma^{1/2} y\|_2\leq \sigma\ones^Ty.
\]
$y=(h,I)$ is the liquid and illiquid exposure. Thus, we use the constraint 
\[
\|\Sigma^{1/2} y\|_2\leq \sigma\ones^Ty,
\]
which is invariant in wealth.  It is also convex, which means that problems with such constraints
can be reliably solved.

\paragraph{Insolvency constraint.}
An important challenge in performing strategic asset allocation with illiquid alternatives
is ensuring that the probability of being unable to pay a capital call is extremely low.
In our model, this corresponds to requiring
\[
P(W_{t+1}<0 \mid X_t,n_t,h_t) \leq \epsilon^\text{ins}
\]
for a small probability of failure $\epsilon^\text{ins}$.  We make several approximations to
facilitate a convex constraint.
First, we approximate $R_t^\text{liq}$ as a multivariate normal random variable,
\[
R_t^\text{liq}\sim N(\mu_\text{liq},\Sigma_\text{liq}).
\] 
It is important to note that these parameters are the mean and covariance of the liquid returns,
rather than the mean and covariance which parameterize the log normal liquid return distribution given by
$\mu_\text{ret}$ and $\Sigma_\text{ret}$ in \eqref{e-join-latent}.
Then we assume we receive the expected calls $\overline{c}_t =\Expect[C_t\mid X_t,n_t,h_t]$, 
which is a linear function of our controls. 
They are given by
  $\overline{c}_t =\overline{\lambda}_t^{1, T}K_t+\overline{\lambda}_t^{0,T}n_t$.
Finally, we assume pessimistically there are no distributions or outside cash.
With these approximations, we have
\begin{eqnarray*}
P(W_{t+1}<0 \mid X_t,n_t,h_t) & \approx &  P(R_t^\text{liq}h_t-\overline{c}_t \leq 0)\\
&=& P(N(h_t^T\mu_\text{liq}-\overline{c}_t,h_t^T\Sigma_\text{liq}h_t)<0) \leq \epsilon.
\end{eqnarray*}
This probabilistic constraint holds if and only if
\begin{equation}\label{e-insolvency}
  \overline{c}_t-h_t^T\mu_\text{liq} \leq \Phi^{-1}(\epsilon^\text{ins})\|\Sigma_\text{liq}^{1/2}h_t\|_2,    
\end{equation}
where $\Phi$ is the standard normal cumulative distribution function. 
This constraint is convex provided $\epsilon^\text{ins}\leq 1/2$,
since then $\Phi^{-1}(\epsilon^\text{ins})\leq 0$, and \eqref{e-insolvency} 
is a second order cone constraint (see \cite[\S 4.4.2]{cvxbook}).

As mentioned above, the constraint \eqref{e-insolvency} is pessimistic because it assumes no distribution. 
An alternative and less pessimistic formulation of the insolvency constraint would consider the 
distribution, the calls, and the liquid returns all under a joint normal approximation. 

\paragraph{Smoothing penalty.}
Among control sequences with similar objective values, we would like for new commitments to be
fairly smooth across time. 
We can consider a natural commitment smoothing penalty
\[
g(n) = \sum_{t=0}^{H-1}\gamma^t\|n_{t+1}-n_t\|^2
\]
The time discount $\gamma$ appears because in a growth context we expect $n_t$ to increase over time.  
Additionally, it helps account for the increased
uncertainty of future planned steps.


\paragraph{MPC planning problem.}
All objective terms and constraints outlined above are consolidated into one optimization problem. 
At time $t$, we plan $\{\hat x_{\tau|t}\}_{\tau=t}^{t+H+1}$, $\{\hat u_{\tau|t}\}_{\tau=t}^{t+H}$, 
where $H$ is the planning horizon, by solving the optimization problem 
\begin{equation}\label{e-full-mpc-problem}
\begin{array}{lll}
  \text{maximize} & \sum_{\tau=t}^{t+H}\gamma^t \paren{\hat L_{\tau|t}
  +\ones^T\hat I_{\tau|t}-\lambda^\text{cash}\hat s_{\tau|t}}-\lambda^\text{smooth}g(\hat n_{\cdot|t})\\
  \text{subject to} &  \hat x_{t|t} = x_t\\
  & \hat x_{\tau+1|t}=\overline{A}
 \hat x_{\tau|t}+
  \overline{B} \hat u_{\tau|t}, & \tau=t,\ldots,t+H\\
  & \hat L_{\tau|t}\geq 0, & \tau= t,\ldots,t+H+1\\
  & \hat h_{\tau|t},\hat n_{\tau|t},\hat s_{\tau|t} \geq 0, & \tau= t,\ldots,t+H\\
  & \mathbf{1}^T \hat h_{\tau|t} = \hat L_{\tau|t}, & \tau=t,\ldots,t+H\\
  & \|\Sigma^{1/2} \hat y_{\tau|t}\|_2\leq \sigma\ones^T\hat y_{\tau|t}, & \tau=t,\ldots,t+H\\
  & \overline{\lambda}^{1,T} \hat K_\tau+\overline{\lambda}^{0,T}\hat n_{\tau|t}
  -\hat h_{\tau|t}^T\mu_\text{liq} 
  \leq \Phi^{-1}(\epsilon^\text{ins})\|\Sigma_\text{liq}^{1/2} \hat h_{\tau|t}\|_2,
  & \tau=t,\ldots,t+H,\\
\end{array}
\end{equation}
where $\lambda^\text{cash}>0$ is a hyperparameter penalizing outside cash use.
Recall that $L$, $I$, and $K$ are components of $x$, and $h$, $n$,
and $s$ are components of $u$.

\subsection{Example}\label{results}    

In this example, we evaluate the performance of the two policies described in
\S\ref{heuristic} and \S\ref{MPC} using the risk return trade off. For the
parameters defined in \eqref{e-join-latent}, we use the specific values of
{
\begin{eqnarray}\label{e-example-ret-liq}
\mu_\text{ret} &=&\pmat{
0.158,&
0.000,&
0.072,&
0.023,&
0.036,&
0.046
},\\
\sigma_\text{ret} &=& \pmat{
0.281,&
0.000,&
0.206,&
0.046,&
0.047,&
0.162
},\\
C_\text{ret} &=& \bmat{
1.000  & 0.000  & 0.422  & -0.298 & -0.002 & 0.261  \\
0.000  & 1.000  & 0.000  & 0.000  & 0.000  & 0.000  \\
0.422  & 0.000  & 1.000  & -0.843 & 0.197  & 0.800  \\
-0.298 & 0.000  & -0.843 & 1.000  & -0.018 & -0.739 \\
-0.002 & 0.000  & 0.197  & -0.018 & 1.000  & 0.628  \\
0.261  & 0.000  & 0.800  & -0.739 & 0.628  & 1.000  
}\label{e-example-cov-liq}
\end{eqnarray}
with $\Sigma_\text{ret} = \diag(\sigma_\text{ret})C\diag(\sigma_\text{ret})$.

\paragraph{Illiquid dynamics.}
We consider the actual illiquid world: full call/distribution random dynamics as
described in \S\ref{stochastic_linear_system_joint}. We evaluate the two policies
described in \S\ref{heuristic} and \S\ref{MPC} on the same simulated
returns as the imaginary Markowitz portfolio.

\paragraph{Example policy specifications.} In this case study, 
we use the steady state commitment policy with parameters 
\[
c=3.685,\qquad \kappa=.1,
\] and values of $\theta$ arising from solving the one period Markowitz problem
defined in \ref{markowitz} for 30 evenly spaced values of $\sigma$ between 0 and
.3, with our specified return distribution parameters defined in
(\ref{e-example-mean-return}--\ref{e-example-cov-return}).

For the MPC policy, we use the same $\sigma$ values described above, but for
numerical reasons use the standard trick of moving the risk limit to penalized
form by subtracting
\[
\lambda^\text{risk}(\|\Sigma^{1/2} y_t\|_2-\sigma\ones^Ty_t)_+
\]
from each term of the objective defined in \eqref{e-full-mpc-problem}, penalizing excess risk.
The parameter values are
\[
\gamma=.97,\qquad
H = 10,\qquad
\epsilon^\text{ins}=.02,\qquad
\lambda^\text{risk}=10,\qquad
\lambda^\text{smooth}=.1,\qquad
\lambda^\text{cash}=1000,
\]
with $\overline{A},\overline{B}$ as defined in \eqref{e-mean-liquid-dynamics}, 
with the distributions instantiated in \eqref{e-intensity-dist}
and (\ref{e-example-ret-liq}--\ref{e-example-cov-liq}).

\paragraph{Results.} We see in figure~\ref{f-policy-results-20} that both the MPC and
heuristic policies are extremely close to the risk-return performance of the
liquid relaxation, which is an unattainable benchmark. This
is despite the challenging illiquid dynamics we face in the non-relaxed setting.
\begin{figure}
  \begin{center}
  \includegraphics[scale=.3]{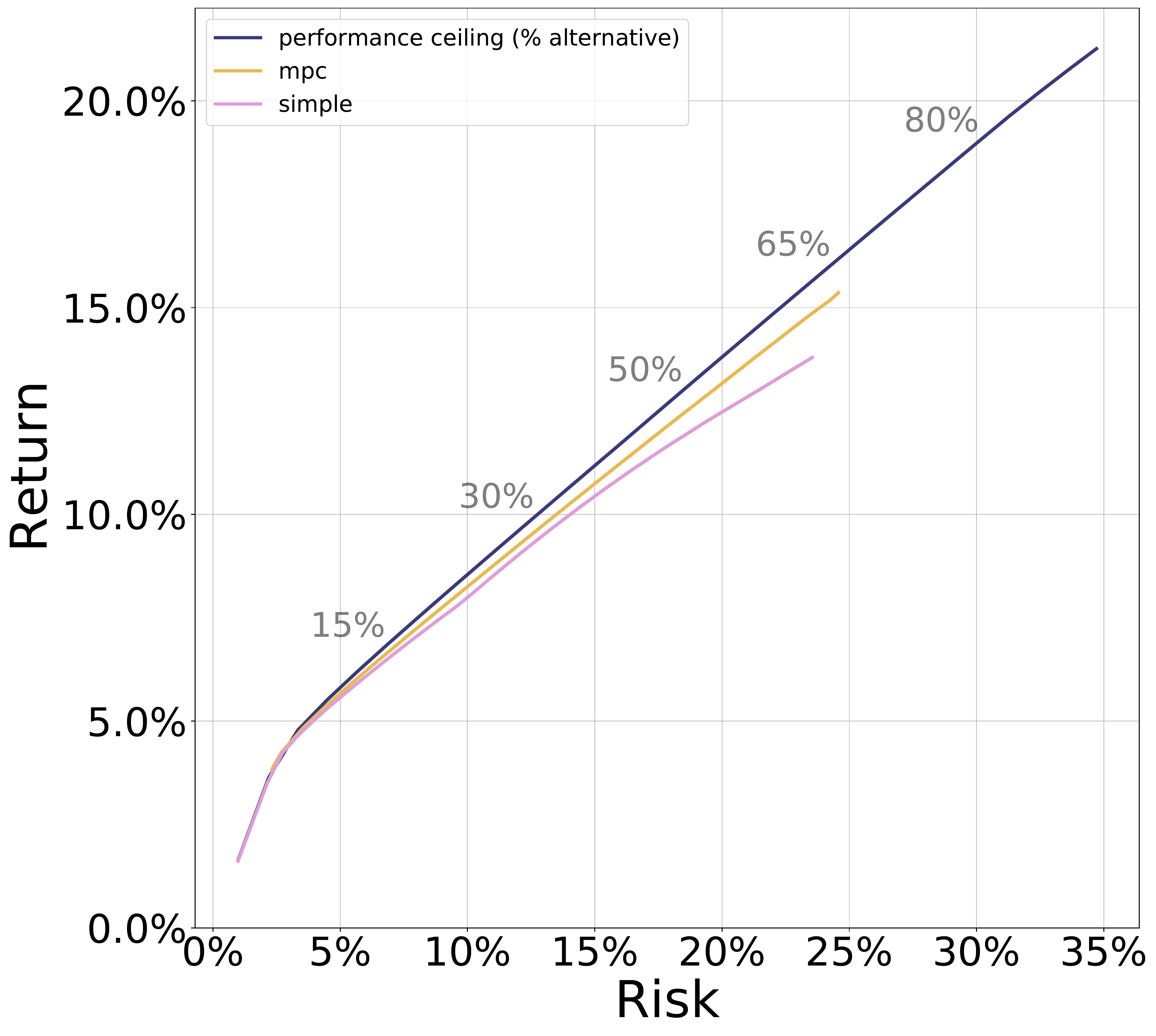}
  \caption{Risk return trade-off, 200 simulations of 20 periods.}
  \label{f-policy-results-20}
  \end{center}
  \end{figure}
The performance stated here is averaged across 20 periods of simulation, 
for 200 simulated trajectories. 

We can also examine the performance across a shorter time horizon. 
figure \ref{f-policy-results-10} shows the same risk-return tradeoff for
10 periods. 
Evidently, there is a larger gap between the MPC policy and the liquid performance
ceiling, and also between the MPC and simple policies. 
This has a perfectly clear interpretation: because there is a roughly 4 period 
delay before peak illiquid exposure (see figure \ref{impulse_response}), the impact 
of the illiquid alternative asset's high returns is delayed. 
Additionally, by planning ahead, the MPC policy achieves illiquid exposure faster
than the simple policy.
\begin{figure}
  \begin{center}
  \includegraphics[scale=.3]{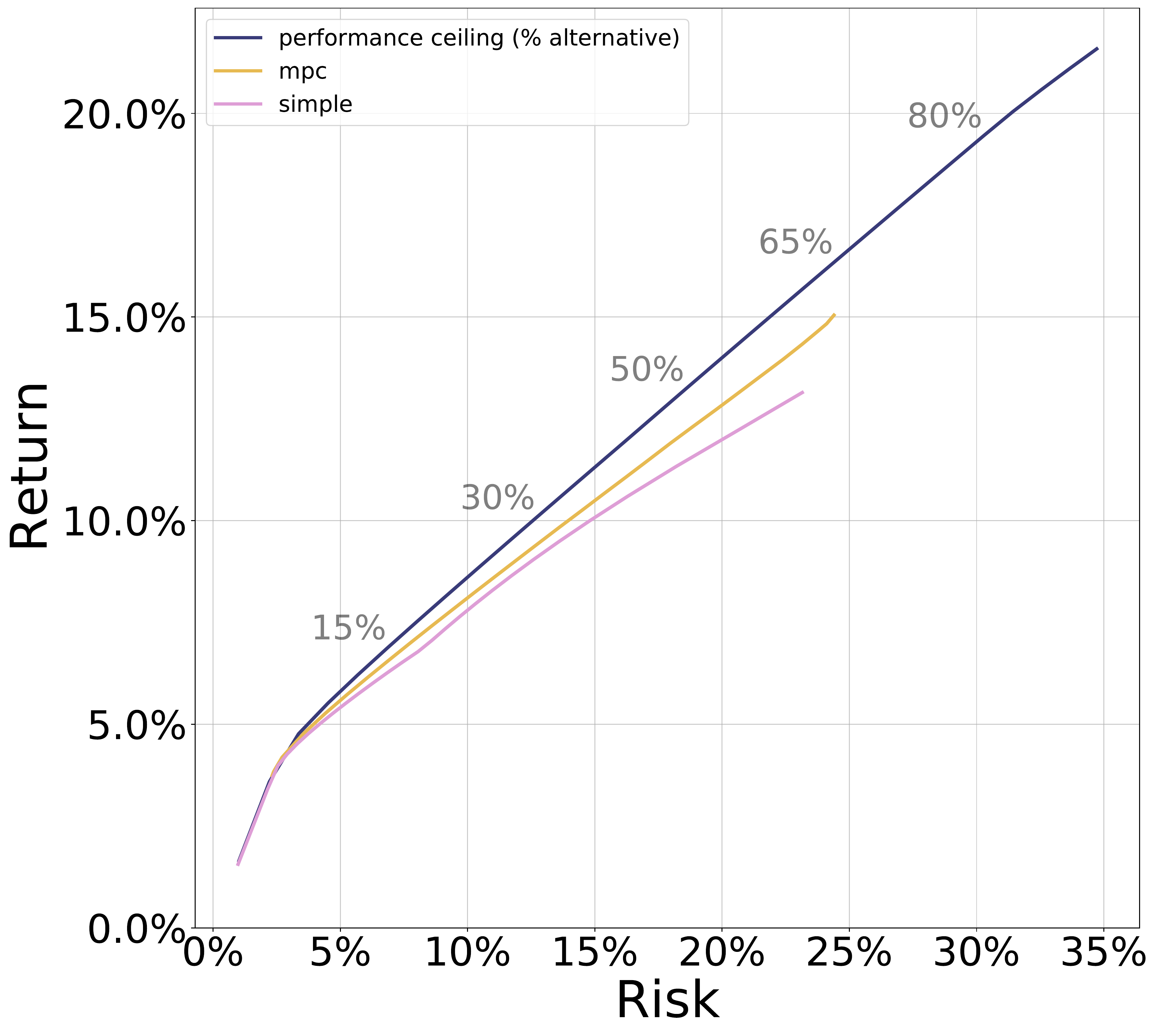}
  \caption{Risk return trade-off, 200 simulations of 10 periods.}
  \label{f-policy-results-10}
  \end{center}
\end{figure}

By looking at the average allocation across time for both policies shown in
figure~\ref{allocations}, we can further understand these differences. 
We can now see directly that the MPC policy is able to reach a stable allocation in fewer
periods than the heuristic policy. If we include the proportional feedback
control, the heuristic does reach the allocation faster, but still not as
quickly as the MPC. Another difference is that the heuristic policy and MPC sweep out the
same risk return trade-off profile, but may not choose the exact same portfolio
steady state weights. Generally, the heuristic undershoots the illiquid target
it is trying to reach.
\begin{figure}
\centering
\includegraphics[scale=.5]{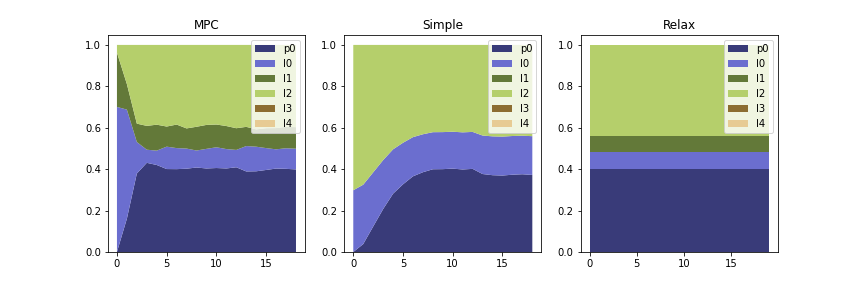}
\caption{Average allocations across time, 200 simulations for 20 periods. 
MPC, heuristic, and relaxation.}
\label{allocations}
\end{figure}

\clearpage
\section{Extensions}\label{s-extensions}
\paragraph{Liquidation.}
We can easily extend the model to allow for liquidation of illiquid alternatives
on the secondaries market. Per common practice (for example, see
\cite{Giommetti2021}), we assume that at time $t$ we can liquidate $0\leq
\ell_t\leq P_t$ from $P_t$ which, after a haircut $\phi$ is available as liquid
wealth $\phi\ell_t$. This changes the control by appending an $\ell_t\in\reals^{n^\text{ill}}$ to $u_t$ as
defined in \S\ref{stochastic_linear_system_joint}. Accordingly, the new control matrix is given by 
\[\bmat{
B_t^{\text{liq}} & \tilde{B}_t
},
\]
with
\[
\tilde{B}_t = \bmat{
\phi \ones^T\\
-I\\ 0
},
\]
where the block of zeros and the identity matrix are in dimension $\reals^{n^\text{ill}\times n^\text{ill}}$.
There is also a new constraint enforcing $0 \leq-\ell_t\leq I_t$, \ie the maximum
liquidation is the entire liquid exposure in a given asset. 

\paragraph{Tracking fixed weights.}
In this paper, a user  specifies a risk tolerance parameter $\sigma$ as in \eqref{e-markowitz-problem},
which implicitly specifies the portfolio weights across the liquid and illiquid assets. However, an investor
may have arrived with pre-selected target portfolio weights. Instead of seeking to track a target illiquid exposure,
as in the problem posed in \S\ref{s-commitment-opt}, we can instead seek to track target portfolio weights.
A natural tracking constraint in planning is  
\[
\|\Sigma^{1/2}(\hat{y}_{\tau|t}-\theta\ones^T\hat{y}_{\tau|t})\|_2
\leq \sigma^{\text{track}}\ones^T\hat{y}_{\tau|t},
\]
where $\theta\in\reals^{n^\text{ill}+n^\text{liq}}$ is the user-provided vector of target portfolio weights, and $\sigma^\text{track}$
is a tracking variance hyperparameter. As with our risk constraint in \eqref{e-full-mpc-problem}, in practice a slack variable can be added
to the above constraint to guarantee feasibility. 

\paragraph{Liabilities.}
We can incorporate external liabilities $Z_t$ by modifying our liquid wealth update to 
\[
L_{t+1}=h_t^T R_t^\text{liq}-\ones^TC_t+\ones^TD_t+s_t-Z_t.
\]
This encodes an external obligation of $Z_t$ dollars in period $t$. 
This could represent the liabilities of an insurer or a pension fund.
MPC is able to handle these liabilities quite gracefully: at every time $t$ the planning problem
takes in a forecast of the next $H$ liabilities $\hat L_{\tau|t}$, $\tau=t,\ldots,t+H$. 
The insolvency constraint~\eqref{e-insolvency} can be modified to include the liabilities as
\[
L_t+\overline{c}_t-h_t^T\mu_\text{liq} \leq \Phi^{-1}(\epsilon^\text{ins})\|\Sigma_\text{liq}^{1/2}h_t\|_2.
\]

\paragraph{Time varying forecasts.}
In the current problem formulation, we plan based on the mean dynamics \eqref{e-mean-liquid-dynamics}, 
which we treat as stationary at every time $t$. The mean dynamics capture the expected returns, call intensities, 
and distribution intensities. It is immediate to replace these stationary forecasts with time varying ones: 
planning at time $t$ in \eqref{e-full-mpc-problem} becomes
\[
\hat x_{\tau+1|t}=\overline{A}_{\tau|t}
\hat x_{\tau|t}+
\overline{B}_{\tau|t} \hat u_{\tau|t},\qquad \tau=t,\ldots,t+H,
\]
where $\overline{A}_{\tau|t}$ and $\overline{B}_{\tau|t}$ are the forecasted mean dynamics at time $\tau$ generated
at time $t$.

\paragraph{Illiquid dynamics with vintages.}
A natural way to extend the Takahashi and Alexander illiquid asset model is to have time varying intensity parameters
that depend on the age of the investment. This amounts to keeping track of vintages for each asset class, rather than 
aggregating all exposure to a given illiquid asset in one state, as this paper does. This extension is quite natural,
and is readily implementable as an only slightly larger tractable convex optimization based planning problem.
A given illiquid asset at time $t$, rather than by two states $I_t$ and $K_t$, will now require $2k$ states, 
\[
I_{t,a},\quad K_{t,a},
\]
where $a$ denotes the age of the investment and the maximum age tracked is $k$. In words, at each time, we keep track of the exposure and uncalled commitments
from commitments of age $a$.

\section{Conclusion}
We have described a flexible stochastic linear system model of liquid and illiquid alternative assets,
that takes into account the dynamics of the illiquid assets and the 
randomness of returns, calls, and distributions.
This model allows us to develop an MPC policy that in each time period chooses 
a liquid wealth allocation, and also new commitments to make in each alternative asset.

We compare the results of this policy with a relaxed liquid model, where we assume
that all illiquid assets are fully liquid.  
This relaxed liquid model is easy to understand, since the challenges of alternative assets
have all been swept under the rug.
For the relaxed liquid model, we can work out optimal investment policies.  The performance
with these policies can be thought of as an unattainable benchmark, 
that we know we cannot achieve or beat when 
all the challenges of alternative investments are present.

Suprisingly, the performance of the MPC policy under the real model, with all the 
challenges of alternative assets, is very close to the performance of the relaxed
liquid model, under an optimal policy.
Roughly speaking, there isn't much room for improvement.
This is a strong validation of the MPC policy.

Another interesting conclusion is that the relaxed liquid model is not as useless as 
one might imagine, since MPC can attain similar performance with all the challenges present.
In a sense this validates reasoning based on the relaxed liquid model, where illiquid assets
are treated as liquid assets.
Roughly speaking, the asset manager can reason about the portfolio using the simple relaxed 
liquid model; feedback control with the MPC policy handles the 
challenges of illiquid alternative assets.

\clearpage
\bibliographystyle{alpha}
\bibliography{library}

\end{document}